\documentclass[11pt]{amsart}
\usepackage{amscd}      
\usepackage{amssymb}
\usepackage{amsmath}
\usepackage{xypic}      
\LaTeXdiagrams          
\usepackage[all,v2]{xy}
\xyoption{2cell} \UseAllTwocells \xyoption{frame} \CompileMatrices
\allowdisplaybreaks[3]

\usepackage{latexsym}
\usepackage{epsfig}
\usepackage{amsfonts}
\usepackage{enumerate}

\newtheorem{prop}{Proposition}[section]
\newtheorem{lem}[prop]{Lemma}

\newtheorem{cor}[prop]{Corollary}
\newtheorem{thm}[prop]{Theorem}

\numberwithin{equation}{section}

\newtheorem{defn}[prop]{Definition}

\newtheorem{example}[prop]{Example}

\newtheorem{rmk}[prop]{Remark}

\newenvironment{pf}{\begin{trivlist}\item[]{\sc Proof.}}%
            {\nolinebreak $\Box$ \end{trivlist}}

\newcommand{\noprint}[1]{}

\newcommand{\toto}{\rightrightarrows}

\newcommand{\GG}{{\mathfrak G}}

\newcommand{\HH}{{\mathfrak H}}

\newcommand{\ldiag}[1]%
       {\makebox[0cm]{${\scriptstyle#1}\downarrow\phantom{\scriptstyle#1}$}}
\newcommand{\ldiagup}[1]%
       {\makebox[0cm]{${\scriptstyle#1}\uparrow\phantom{\scriptstyle#1}$}}
\newcommand{\rdiag}[1]%
       {\makebox[0cm]{$\phantom{\scriptstyle#1}\downarrow{\scriptstyle#1}$}}
\newcommand{\sediagr}[1]%
       {\makebox[0cm]{$\phantom{\scriptstyle#1}\searrow{\scriptstyle#1}$}}
\newcommand{\nediagr}[1]%
       {\makebox[0cm]{$\phantom{\scriptstyle#1}\nearrow{\scriptstyle#1}$}}
\newcommand{\rdiagup}[1]%
       {\makebox[0cm]{$\phantom{\scriptstyle#1}\uparrow{\scriptstyle#1}$}}
\newcommand{\swdiag}[1]%
       {\makebox[0cm]{$\phantom{\scriptstyle#1}\swarrow{\scriptstyle#1}$}}
\newcommand{\sediag}[1]%
       {\makebox[0cm]{${\scriptstyle#1}\searrow\phantom{\scriptstyle#1}$}}
\newcommand{\nediag}[1]%
       {\makebox[0cm]{${\scriptstyle#1}\nearrow\phantom{\scriptstyle#1}$}}

\newcommand{\doublearrowstack}[2]%
                      {{{{\scriptstyle#1}\atop{\textstyle\longrightarrow}}\atop{{\textstyle\longrightarrow}\atop{\scriptstyle#2}}}}
\newcommand{\rightleftarrowstack}[2]%
                      {{{{\scriptstyle#1}\atop{\textstyle\longrightarrow}}\atop{{\textstyle\longleftarrow}\atop{\scriptstyle#2}}}}
\newcommand{\leftrightarrowstack}[2]%
                      {{{{\scriptstyle#1}\atop{\textstyle\longleftarrow}}\atop{{\textstyle\longrightarrow}\atop{\scriptstyle#2}}}}

\newcommand{\overtoparrow}%
{\makebox[0cm]{\beginpicture \setcoordinatesystem units
<.8cm,.4cm> point at 0 0 \setplotarea x from -3 to 3, y from 0 to
1 \setquadratic \plot -3 0 0 1 3 0 / \put{\vector(3,-1){0}}[Bl] at
3 0
\endpicture}}

\newcommand{\underbottomarrow}%
{\makebox[0cm]{\beginpicture \setcoordinatesystem units
<.8cm,.4cm> point at 0 0 \setplotarea x from -3 to 3, y from 0 to
1 \setquadratic \plot -3 1 0 0 3 1 / \put{\vector(3,1){0}}[Bl] at
3 1
\endpicture}}

\newcommand{\ses}[5]%
{0\longrightarrow#1\stackrel{#2}{ \longrightarrow}#3\stackrel{#4}{
\longrightarrow}#5\longrightarrow0}

\newcommand{\dt}[6]%
{#1\stackrel{#2}{longrightarrow}#3
\stackrel{#4}{\longrightarrow}#5 \stackrel{#6}{\longrightarrow}
#1[1]}

\newcommand{\cat}[1]%
{(\mbox{\rm #1})}

\def\Label#1{\label{#1}{\tt [#1]}\phantom{h}}
\def\Label{\label}

\title[Automorphism group of toric DM stacks]{The Automorphism Group of Toric Deligne-Mumford Stacks}

\author{Yunfeng Jiang}
\address{Department of Mathematics\\ University of British Columbia\\ 1984 Mathematics Road\\Vancouver\\ BC V6T 1Z2\\ Canada}
\email{jiangyf@math.ubc.ca}
\address{$Current~ Address$: Department of Mathematics\\ University of Utah\\ 155 S 1400 E JWB 233\\Salt Lake city\\ UT 84112\\ USA}

\date{\today}

\begin{document}
\begin{abstract}
We prove that the automorphism group of a toric Deligne-Mumford stack
is isomorphic to the $2$-group associated to the stacky fan.
\end{abstract}

\maketitle

\section{Introduction}

Let $\Sigma$ be a simplicial fan with $n$ rays in a lattice $N$ of
dimension $d$. The simplicial toric variety $X(\Sigma)$ can be
represented as a geometric quotient $Z/G$, where
$Z:=\mathbb{C}^{n}\setminus V(J_{\Sigma})$, $V(J_{\Sigma})$ is the
closed subvariety defined by the irrelevant ideal $J_{\Sigma}$ of
the fan $\Sigma$, see \cite{Cox}. The group
$G=\mbox{Hom}(A_{d-1}(X(\Sigma)),\mathbb{C}^{*})$ acts on $Z$
through the map $\alpha: G\rightarrow (\mathbb{C}^{*})^{n}$ obtained
by taking $\mbox{Hom}(-,\mathbb{C}^{*})$-functor on the map
$\beta^{\vee}$ in the following exact sequence
\begin{equation}\Label{exact1}
0\longrightarrow N^{\star}\longrightarrow
\mathbb{Z}^{n}\stackrel{\beta^{\vee}}{\longrightarrow}
A_{d-1}(X(\Sigma))\longrightarrow 0,
\end{equation}
in \cite{F}, where $A_{d-1}(X(\Sigma))$ is the Weil divisor class
group of the toric variety $X(\Sigma)$. Let $\mbox{Aut}(X(\Sigma))$
be the automorphism group of $X(\Sigma)$ and $\mbox{Aut}(Z)$ the
automorphism group of $Z$. Since the group $G$ acts on $Z$, it is
naturally a subgroup of the automorphism group $\mbox{Aut}(Z)$. Let
$\widetilde{\mbox{Aut}(Z)}$ be the normalizer of $G$ in
$\mbox{Aut}(Z)$. Then from \cite{Cox} we have the following exact
sequence
\begin{equation}\Label{exact2}
1\longrightarrow G\longrightarrow
\widetilde{\mbox{Aut}(Z)}\stackrel{}{\longrightarrow}
\mbox{Aut}(X(\Sigma))\longrightarrow 1.
\end{equation}
If $\widetilde{\mbox{Aut}^{0}(Z)}$ and $\mbox{Aut}^{0}(X(\Sigma))$
are the components of the identity elements in
$\widetilde{\mbox{Aut}(Z)}$ and $\mbox{Aut}(X(\Sigma))$, then we
have the exact sequence
\begin{equation}\Label{exact3}
1\longrightarrow G\longrightarrow
\widetilde{\mbox{Aut}^{0}(Z)}\stackrel{}{\longrightarrow}
\mbox{Aut}^{0}(X(\Sigma))\longrightarrow 1.
\end{equation}

The main goal of this note is to generalize the above results to
toric Deligne-Mumfords stacks. Generalizing the idea of Cox on
simplicial toric varieties, Borisov, Chen and Smith \cite{BCS}
defined the notion of toric Deligne-Mumford stacks which are encoded
in terms of  stacky fans. A stacky fan
$\mathbf{\Sigma}=(N,\Sigma,\beta)$ is a triple, where $N$ is a
finitely generated abelian group, $\Sigma\subset
N\otimes_{\mathbb{Z}}\mathbb{Q}$ is a simplicial fan and $\beta:
\mathbb{Z}^{n}\longrightarrow N$ is a map determined by the elements
$\{b_{1},\cdots,b_{n}\}$ in $N$.  We require that $\beta$ has finite
cokernel and $\{\overline{b}_{1},\cdots,\overline{b}_{n}\}$ generate
the simplicial fan $\Sigma$, where $\overline{b}_{i}$ is the image
of $b_{i}$ under the natural map $N\longrightarrow
\overline{N}=N/N_{tor}$. Then we have the exact sequences:
\begin{equation}\Label{exact4}
0\longrightarrow DG(\beta)^{\star}\longrightarrow
\mathbb{Z}^{n}\stackrel{\beta}{\longrightarrow} N\longrightarrow
Coker(\beta)\longrightarrow 0,
\end{equation}
and
\begin{equation}\Label{exact5}
0\longrightarrow N^{\star}\longrightarrow
\mathbb{Z}^{n}\stackrel{\beta^{\vee}}{\longrightarrow}
DG(\beta)\longrightarrow Coker(\beta^{\vee})\longrightarrow 0,
\end{equation}
where $\beta^{\vee}$ is the Gale dual of $\beta$ (see \cite{BCS}).
The  toric Deligne-Mumford stack $\mathcal{X}(\mathbf{\Sigma})$
associated to $\mathbf{\Sigma}$ is defined to be  the quotient
stack $[Z/G]$, where $Z$ is the same as in the quotient
construction of toric varieties,
$G=\mbox{Hom}(DG(\beta),\mathbb{C}^{*})$ and the action is through
a group homomorphism $\alpha: G \longrightarrow
(\mathbb{C}^{*})^{n}$ in the exact sequence
\begin{equation}\Label{exact6}
1\longrightarrow \mu\longrightarrow
G\stackrel{\alpha}{\longrightarrow}
(\mathbb{C}^{*})^{n}\longrightarrow
(\mathbb{C}^{*})^{d}\longrightarrow 1, \
\end{equation}
which is obtained by taking $\mbox{Hom}(-,\mathbb{C}^{*})$-functor
on the exact sequence (\ref{exact5}).

Let $\mathcal{X}$ be an algebraic stack. The automorphism group of
$\mathcal{X}$ is naturally a $2$-group because an isomorphism of the
stack is a $2$-isomorphism. A 2-group is a group object in the
category of groupoids. The notion of $2$-groups, also called
crossed-modules, appeared first in algebraic topology. The
fundamental concepts of $2$-groups and crossed-models can be found
in \cite{Noohi}. A crossed-module $\GG=[G_2\rightarrow G_1]$
consists of  a pair of groups $G_1$ and $G_2$, a group homomorphism
$\varphi: G_2\rightarrow G_1$, and a (right) action of $G_1$ on
$G_2$, denoted by $^{-}a$, which lifts the conjugation action of
$G_1$ on the image of $\varphi$ and descends the conjugation action
of $G_2$ on itself.

Let $\mathcal{X}(\mathbf{\Sigma})$ be  a toric Deligne-Mumford
stack associated to the stacky fan $\mathbf{\Sigma}$. Since
$(\mathbb{C}^{*})^{n}$ is the maximal torus of the automorphism
group $\mbox{Aut}(Z)$,  the map $\alpha$ in the exact sequence
(\ref{exact6}) naturally induces a map $\varphi: G\rightarrow
\widetilde{\mbox{Aut}(Z)}$ whose image lies in the centralizer
$\widetilde{\mbox{Aut}^{0}(Z)}$ of $G$ in
$\widetilde{\mbox{Aut}(Z)}$. We define the action of
$\widetilde{\mbox{Aut}(Z)}$ on $G$ by:
\begin{equation}\Label{action}
\begin{cases}
h\cdot g=g&\text{if
}~ h\in \widetilde{\mbox{Aut}^{0}(Z)};\\ \\
h\cdot g=hgh^{-1}&\text{if}~h\in \widetilde{\mbox{Aut}(Z)}\setminus
\widetilde{\mbox{Aut}^{0}(Z)}\,.
\end{cases}
\end{equation}
Then
$$\GG=[G\stackrel{\varphi}{\rightarrow} \widetilde{\mbox{Aut}(Z)}]$$
is a crossed module and it defines a (weak) $2$-group. We call
$\GG$ the 2-group associated to the stacky fan $\mathbf{\Sigma}$.
Let $\mbox{Aut}(\mathcal{X}(\mathbf{\Sigma}))$ be the (weak)
$2$-group of the automorphism group of the toric Deligne-Mumford
stack $\mathcal{X}(\mathbf{\Sigma})$. Then there is a natural
homomorphism from $\GG$ to
$\mbox{Aut}(\mathcal{X}(\mathbf{\Sigma}))$ which will be discussed
in Section 2.

\begin{thm}
Let $\mathcal{X}(\mathbf{\Sigma})$ be  a toric Deligne-Mumford
stack associated to the stacky fan $\mathbf{\Sigma}$. Then the
natural map
$$f: \GG\rightarrow \mbox{Aut}(\mathcal{X}(\mathbf{\Sigma}))$$
is an isomorphism.
\end{thm}

Actually we prove a result for general quotient stack $[X/G]$ in
which the group $G$ is abelian such that the theorem is a corollary
of this result. Our result generalizes Lemma 8.2 of Behrend and
Noohi \cite{BN},where they only consider the case of the centralizer
of the group $G$ in the automorphism group of $X$.

The paper is outlined as follows. In Section \ref{2group} we quickly
review the notion of $2$-groups.  The main Theorem  is proved in
Section \ref{theorem}. In Section \ref{weighted} we prove that the
weighted projective stacks are toric Deligne-Mumford stacks and
discuss the weighted projective linear $2$-groups which are the
automorphism group of weighted projective stacks.

\subsection*{Conventions}
In this paper we work entirely algebraically over the field of
complex numbers.  By an orbifold we mean a smooth
Deligne-Mumford stack with trivial generic stabilizer.

We write $N^{\star}$ for $Hom_{\mathbb{Z}}(N,\mathbb{Z})$ and $N\to
\overline{N}$ the natural map of modulo torsion. We refer the reader
to \cite{BCS} for the construction of the Gale dual $\beta^{\vee}:
\mathbb{Z}^{n}\to DG(\beta)$ of $\beta: \mathbb{Z}^{n}\to N$.

We use the notation $\GG$ to represent $2$-groups. For more details about
$2$-groups and crossed modules, see \cite{Noohi}.

\subsection*{Acknowledgments}
I would like to thank my advisor Kai Behrend for encouragements
and valuable discussions.
\section{$2$-groups and automorphism group of stacks}\Label{2group}
In this section we review the notions of $2$-groups and crossed modules.
The automorphism group of an algebraic  stack is
naturally a $2$-group.

\subsection{$2$-groups}\label{2-group}
A 2-group $\GG$ is a group object in the category of groupoids.
Alternatively, we can define a 2-group to be a groupoid object in
the category of groups, or also, as a (strict) 2-category with one
object in which all 1-morphisms and 2-morphisms are invertible (in
the strict sense). If we require that the 1-morphisms are only
equivalences and not necessarily strictly invertible, what we obtain
is called a weak 2-group.

A morphism $f : \GG\rightarrow \HH$ of 2-groups is a map of
groupoids that respects the group operation on the nose. If we view
$\GG$ and $\HH$ as 2-categories with one object, such f is nothing
but a strict 2-functor. Same definition applies to weak 2-groups
(because we have not weakened associativity). To a 2-group $\GG$ we
associate the groups $\pi_1(\GG)$ and $\pi_2(\GG)$ as follows. The
group $\pi_1(\GG)$ is the set of isomorphism classes of object of
the groupoid $\GG$. The group structure on $\pi_1(\GG)$ is induced
from the group operation of $\GG$. The group $\pi_2(\GG)$ is the
group of automorphisms of the identity object $e$ in $\GG$. This is
an abelian group. Any morphism $f : \GG\rightarrow\HH$ of (weak)
2-groups induces homomorphisms on $\pi_1(\GG)$ and $\pi_2(\GG)$. We
say such $f$ is an $equivalence$ if both these maps are
isomorphisms. (Warning: an equivalence need not have an inverse.
Also, two equivalent 2-groups may not be related by an equivalence,
but just a zig-zag of equivalences.) We are usually interested in
2-groups up to equivalence, so we will think of a 2-group as an
object in the homotopy category of 2-groups. This category is
defined by taking the category of 2-groups and inverting the
equivalences. For the readers information, we point out that, there
is a model structure on the category of 2-groupoids; 2-groups are
the pointed connected objects in this category.

\subsection{Crossed moduels}
A crossed-module $\GG = [\varphi: G_2\rightarrow G_1]$ consists of a
pair of groups $G_1$ and $G_2$, a group homomorphism $\varphi:
G_2\rightarrow G_1$, and a (right) action of $G_1$ on $G_2$, denoted
by $^{-}g$, that lifts the conjugation action of $G_1$ on the image
of $\varphi$ and descends the conjugation action of $G_2$ on itself.
The kernel of $\varphi$ is a central (in particular abelian)
subgroup of $G_2$ and is denoted by $\pi_2(\GG)$. The image of
$\varphi$ is a normal subgroup of $G_1$ whose cokernel is denoted by
$\pi_1(\GG)$. A (strict) morphism of crossed-modules is a pair of
group homomorphisms which commute with the $\varphi$ maps and
respect the action of $G_1$ on $G_2$. A morphism is called an
equivalence if it induces an isomorphism on $\pi_1(\GG)$ and
$\pi_2(\GG)$.

\textbf{Equivalence of 2-groups and crossed-modules.} There is a
natural pair of inverse equivalences of categories between the
category $\mathbf{2Gp}$ of 2-groups and the category
$\mathbf{CrossedMod}$ of crossed-modules. Furthermore, these
functors preserve $\pi_1(\GG)$ and $\pi_2(\GG)$. Here is how these
functors are defined.

\textbf{Functor from 2-groups to crossed-modules.} Let $\GG$ be a
2-group. Let $G_1$ be the group of objects of $\GG$, and let $G_2$
be the set of all arrows emanating from the identity object $e$;
$G_2$ is also a group (namely, it is a subgroup of the group of
arrows of $\GG$). Define
$$\varphi: G_2\rightarrow G_1$$
by $\varphi(\alpha):= t(\alpha)$, for $\alpha\in G_2$. The action
of $G_1$ on $G_2$ is given by conjugation. That is, given
$\alpha\in G_2$ and $g\in G_1$, the action is given by
$g^{-1}\alpha g$.

\textbf{Functor from crossed-modules to 2-groups.} Let $[\varphi:
G_2\rightarrow G_1]$ be a crossed-module. Consider the groupoid
$\GG$ whose underlying set of objects is $G_1$ and whose set of
arrows is $G_1\times G_2$. The source and target maps are given by
$s(g, \alpha)=g$, $t(g, \alpha)=g\varphi(\alpha)$. Two arrows $(g,
\alpha)$ and $(h, \beta)$ such that $g\varphi(\alpha)=h$ are
composed to $(g, \alpha\beta)$. The group operation on the set of
objects $\mathbf{ObG}=G_1$ is naturally extended to a group
operation on $\GG$ by setting $(g,\alpha)(h,\beta)=(gh,
(\alpha^{h})\beta)$, where $-^{h}$ stands for the action of $G_1$
on $G_2$.
\subsection{Automorphism group of stacks}
Let $\mathcal{X}$ be a smooth Deligne-Mumford stack. From
\cite{l-mb}, let $X\longrightarrow \mathcal{X}$ be an etale
presentation of the stack $\mathcal{X}$, then $\mathcal{X}$ can be
represented by a etale groupoid
$\mathcal{X}=[X\times_{\mathcal{X}}X\toto X]$ with source and target
maps $s,t$.  For the groupoid $\mathcal{X}$, the automorphism group
$\GG=\mbox{Aut}(\mathcal{X})$ is a (weak) 2-group. We have the
associated groups $\pi_{1}(\GG)$ and $\pi_{2}(\GG)$ which are
defined in Section \ref{2-group}.

We are mainly interested in the quotient stack $\mathcal{X}=[X/G]$
in which the group $G$ is an abelian group. So the associated
groupoid is
$$X\times G\toto X.$$
Let $\mbox{Aut}(X)$ be the automorphism group of $X$, then we have
an embedding $G\longrightarrow \mbox{Aut}(X)$. Let
$\widetilde{\mbox{Aut}(X)}$ be the normalizer of $G$ in
$\mbox{Aut}(X)$ and $\widetilde{\mbox{Aut}^{0}(X)}$  the centralizer
of $G$ in $\mbox{Aut}(X)$, i.e. the component of the identity
element in $\widetilde{\mbox{Aut}(X)}$. We define the action of
$\widetilde{\mbox{Aut}(X)}$ on $G$ as follows: for $g\in G$,
\begin{equation}\Label{action-2.3}
\begin{cases}
h\cdot g=g&\text{if
}~ h\in \widetilde{\mbox{Aut}^{0}(X)};\\ \\
h\cdot g=hgh^{-1}&\text{if}~h\in \widetilde{\mbox{Aut}(X)}\setminus
\widetilde{\mbox{Aut}^{0}(X)}\,;
\end{cases}
\end{equation}
Then we have a 2-group
$$\GG:=[G\rightarrow
\widetilde{\mbox{Aut}(X)}].$$ From \cite{BN}, there is a natural
morphism
$$f: \GG\longrightarrow \mbox{Aut}(\mathcal{X}).$$
In Section \ref{theorem}, we prove that this morphism $\varphi$ is
an isomorphism.

\section{The proof of the Theorem}\label{theorem}
In this section we prove the main theorem of this note. Our proof
bases on a generalization of Lemma 8.2 of Behrend and Noohi
\cite{BN} such that  our main result is a corollary of the
generalization. We first recall the key lemma of Behrend and Noohi
for the automorphism group of abelian quotient stacks.

\begin{lem}\Label{BN} (\cite{BN}) Let $G$ be an abelian group scheme acting on
a connected scheme $X$ over $\mathbb{C}$, and let
$\mathcal{X}=[X/G]$ be the quotient stack. Let
$\widetilde{\mbox{Aut}^{0}(X)}$ be the centralizer of $G$ in the
automorphism group $\mbox{Aut}(X)$ and $\mbox{Aut}(\mathcal{X})$
the (weak) $2$-group of automorphisms of $\mathcal{X}$. Then
\begin{enumerate}
\item The natural homomorphism $\varphi: G\rightarrow \widetilde{\mbox{Aut}^{0}(X)}$
can be turned into a crossed moduel by taking the trivial action of
$\widetilde{\mbox{Aut}^{0}(X)}$ on $G$;
\item If $\GG=[G\stackrel{\varphi}{\longrightarrow}\widetilde{\mbox{Aut}^{0}(X)}]$
is the $2$-group associated to the crossed moduel. Then there is a
natural map of $2$-groups $\GG\rightarrow \mbox{Aut}(\mathcal{X})$.
Furthermore, this morphism induces an isomorphism on $\pi_{2}$;
\item Assume that $\mathcal{X}$ is a proper Deligne-Mumford stack and
$G$ is affine. Then the induced map on $\pi_{1}$ is injective. $\square$
\end{enumerate}
\end{lem}

We construct a new crossed moduel from the quotient stack so that
the associated $2$-group is equivalent to the (weak) automorphism group of the
quotient stack.

\begin{lem}\Label{jiang} Let $G$ be an abelian group scheme acting on
a connected scheme $X$ over $\mathbb{C}$, and let
$\mathcal{X}=[X/G]$ be the quotient stack. Let
$\widetilde{\mbox{Aut}(X)}$ be the normalizer of $G$ in the
automorphism group $\mbox{Aut}(X)$ and  $\mbox{Aut}(\mathcal{X})$
the (weak) $2$-group of automorphisms of $\mathcal{X}$. Then
\begin{enumerate}
\item The natural homomorphism $\varphi: G\rightarrow \widetilde{\mbox{Aut}(X)}$
can be turned into a crossed moduel by taking the action
(\ref{action-2.3}) of $\widetilde{\mbox{Aut}(X)}$ on $G$.
\item If $\GG=[G\stackrel{\varphi}{\longrightarrow}\widetilde{\mbox{Aut}(X)}]$
is the $2$-group associated to the crossed moduel. Then the natural map
of $2$-groups $f: \GG\rightarrow \mbox{Aut}(\mathcal{X})$  is an equivalence
as  $2$-groups.
\end{enumerate}
\end{lem}
\begin{pf}
\begin{enumerate}
\item Since $\varphi$ maps $G$ to the center and the given action of
$\widetilde{\mbox{Aut}(X)}$ on $G$ extends the action of $G$ on itself by conjugation to
$\widetilde{\mbox{Aut}(X)}$, from the definition of crossed moduel,
$\GG$ is a crossed moduel which corresponds to a $2$-group;
\item To prove that the natural map $f: \GG\rightarrow \mbox{Aut}(\mathcal{X})$
is an isomorphism, it suffices to prove that the map induces isomorphisms
on $\pi_2$ and $\pi_1$. We first consider the $\pi_2$ case. We prove that any element
$g\in \mbox{Aut}(\mathcal{X})$ induces an isomorphism of $\mathcal{X}$. Recall that
$\mathcal{X}=[X/G]$ is a quotient stack,
$$ob[X/G](S)=\{(T,\alpha)|T~ \mbox{a}~ G~\mbox{torsor over}~S, \alpha:T\rightarrow
X~\mbox{an}~ G~ \mbox{map}\},$$ and
$$Mor[X/G](S)((T,\alpha),(T^{'},\alpha^{'}))=\{f: T\rightarrow T^{'}
~\mbox{an}~G~\mbox{torsor map such that}~\alpha^{'}\circ
f=\alpha\}.$$ If $g\in \mbox{Aut}^{0}(\mathcal{X})$, then the
induced automorphism of $\mathcal{X}$ is given by keeping the same
torsor $T$ and compose $\alpha$ with the action of $g$ on $X$
which is the same as in \cite{BN}. If $g\in
\mbox{Aut}(\mathcal{X})\setminus \mbox{Aut}^{0}(\mathcal{X})$, the
we have the diagram:
\begin{equation}\Label{map1}
\vcenter{\xymatrix{
T\dto_{\alpha}\rto^{f} & \ T^{'}\dlto^{g\circ\alpha}\\
X, &}}
\end{equation}
then $g\circ\alpha$ is a new $G$-torsor
$T^{'}\stackrel{\alpha^{'}}{\longrightarrow} X$, which we have an
isomorphism between $T$ to $T^{'}$. Hence we can define a morphism
of $2$-groups
$$f: \GG\rightarrow \mbox{Aut}(\mathcal{X}).$$
The proof of isomorphism on $\pi_{2}$ is the same as in
\cite{BN}, we omit the details.

We now prove the isomorphism on $\pi_{1}$. The injectivity on $\pi_{1}$ is the same as in
\cite{BN}. We prove the surjectivity. Let $g$ be an automorphism of the isomorphism classes
of the objects on $\mathcal{X}$, i.e. we have a diagram:
\begin{equation}\label{map2}
\vcenter{\xymatrix{
T\dto_{\alpha}\rto^{\cong} & \ T^{'}\dlto^{g\circ\alpha}\\
X. &}}
\end{equation}
If $T$ is a $G$-torsor over $S$, and
$\alpha: T\rightarrow X$ is  the $G$-map, then $g\circ\alpha$ is also
a $G$-torsor over $S$ and there is an isomorphism between $T$ and $T^{'}$.
So $g\in\widetilde{\mbox{Aut}(X)}$ and the induced map on $\pi_{1}$
is also an isomorphism. So the natural map
$f$ is an equivalence.
\end{enumerate}~~
\end{pf}

\subsection*{The proof of Theorem 1.1}
Let $\mathcal{X}(\mathbf{\Sigma})$ be a toric Deligne-Mumford
stack associated to the stacky fan $\mathbf{\Sigma}$. It is a
quotient stack $[Z/G]$, where $G$ is an abelian group scheme and
affine. So from Lemma \ref{jiang}, the natural map
$$f: \GG\rightarrow \mbox{Aut}(\mathcal{X}(\mathbf{\Sigma}))$$
in the introduction is an equivalence as $2$-groups. $\square$

\section{The weighted projective linear $2$-group}\label{weighted}
In this section we consider the weighted projective stacks. Their
automorphism groups are called weighted projective linear
2-groups.

\subsection{The weighted projective stacks}\label{weightedstack}

\begin{defn} Let $Q=(q_{0},\cdots,q_{n})$ be
a $(n+1)$-tuple of positive integers. The weighted projective
stack of type $Q$,
$\mbox{\bf{P}}(Q)=\mbox{\bf{P}}^{n}_{q_{0},\cdots,q_{n}}$ is
defined to be the quotient stack
$${\bf P}(Q):=[\mathbb{C}^{n+1}\setminus \{0\}/\mathbb{C}^{*}],$$
where the action is given by $\lambda
(x_0,\cdots,x_n)=(\lambda^{q_{0}}x_0,\cdots,\lambda^{q_{n}}x_n)$.
\end{defn}

\begin{rmk}
(1) The above $\mathbb{C}^{*}$-action is free iff $q_{i}=1$ for
every $i=0,\cdots,n$. In this case this is the projective  space
$\mathbf{P}^{n}$. (2) If $gcd(q_{0},\cdots,q_{n})=d\neq 1$, then
$\mbox{\bf{P}}^{n}_{q_{0},\cdots,q_{n}}$ is $reduced$. Otherwise
there exists a gerbe structure over the underlying orbifold.
\end{rmk}

Next we construct weighted projective stacks as toric
Deligne-Mumford stacks. Let $d=gcd(q_{0},\cdots,q_{n})$,
$a_0=\frac{q_0}{d},\cdots,a_n=\frac{q_n}{d}$ and
$$Q_{red}=(a_0,\cdots,a_n).$$
Then $\mathbf{P}(Q_{red})$ is a simplicial toric orbifold. And we
can construct a fan $\Sigma$ in a rank $n$ lattice
$\mathbb{Z}^{n}$. Actually from Fulton \cite{F}, we construct the
simplicial fan $\Sigma$ of $\mathbf{P}(Q_{red})$ as follows. Let
the fan $\Sigma$ be generated by vectors $\{v_{0},\cdots,v_{n}\}$
so that $a_{0}v_{0}+a_{1}v_{1}+\cdots+a_{n}v_{n}=0$, then the
toric variety $X_{\Sigma}$ is the coarse moduli space of the
weighted projective stack $\mathbf{P}(Q_{red})$.

\begin{rmk} The paper of Conrads \cite{Con} gives a method to compute the lattice vectors
$\{v_{0},\cdots,v_{n}\}$ such that they generate a simplicial fan
$\Sigma$ for the corresponding weighted projective stack
$\mathbf{P}(Q_{red})$. Then
$\mathbf{\Sigma_{red}}=(\mathbb{Z}^{n},\Sigma,\overline{\beta})$,
where $\overline{\beta}:\mathbb{Z}^{n+1}\longrightarrow
\mathbb{Z}^{n}$ is given by $\{v_{0},\cdots,v_{n}\}$, is a stacky
fan and $$\mathcal{X}(\mathbf{\Sigma_{red}})=\mathbf{P}(Q_{red}).$$
\end{rmk}

Let $N=\mathbb{Z}^{n}\oplus\mathbb{Z}_{d}$ and consider the map
\begin{equation}\Label{map1}
\beta: \mathbb{Z}^{n+1}\rightarrow N
\end{equation}
given by vectors $\{b_0=(v_{0},1),b_1=(v_1,0),\cdots,b_n=(v_n,0)\}\subset N$.
Then
$$\mathbf{\Sigma}=(N,\Sigma,\beta)$$
is a stacky fan in the sense of Borisov- Chen-Smith \cite{BCS}. From
the definition  of Gale dual in \cite{BCS} we get the following two
diagrams:
\begin{equation}\label{diagram-section4-1}
\begin{CD}
0 @ >>>\mathbb{Z}@ >>> \mathbb{Z}^{n+1}@
>{\beta}>> N @
>>>0\\
&& @VV{}V@VV{}V@VV{}V \\
0@ >>> \mathbb{Z} @ >{}>>\mathbb{Z}^{n+1}@
>{\overline{\beta}}>> \mathbb{Z}^{n}@>>>0,
\end{CD}
\end{equation}
and
\begin{equation}\label{diagram-section4}
\begin{CD}
0 @ >>>\mathbb{Z}^{n}@ >>> \mathbb{Z}^{n+1}@
>{\overline{\beta}^{\vee}}>> \mathbb{Z} @
>>>0@>>> 0\\
&& @VV{}V@VV{}V@VV{[d]}V@V{}VV \\
0@ >>> \mathbb{Z}^{n}@ >{}>>\mathbb{Z}^{n+1}@ >{\beta^{\vee}}>>
\mathbb{Z}@>>>\mathbb{Z}_{d}@>>>0,
\end{CD}
\end{equation}
where $\overline{\beta}^{\vee}$ is given by the matrix
$$[a_0,\cdots,a_{n}]$$
and $\beta^{\vee}$ is given by the matrix
$$[q_{0},\cdots,q_{n}].$$
Applying $\mbox{Hom}_\mathbb{Z}(-,\mathbb{C}^{*})$ to
(\ref{diagram-section4}) yields
\begin{equation}\Label{exactweighted3}
1\longrightarrow \mu_{d}\longrightarrow
\mathbb{C}^{*}\stackrel{\alpha}{\longrightarrow}
(\mathbb{C}^{*})^{n+1}\longrightarrow
(\mathbb{C}^{*})^{n}\longrightarrow 1,\
\end{equation}
where $\alpha$ is given by the matrix
$$\left[\begin{array}{c}
\lambda^{q_{0}}\\
\vdots\\
\lambda^{q_{n}} \end{array} \right].
$$
From the construction of toric Deligne-Mumford stack in the
introduction the weighted projective stack $\mathbf{P}(Q)$ is the
toric Deligne-Mumford stack
$$\mathcal{X}(\mathbf{\Sigma})=[\mathbb{C}^{n+1}\setminus
\{0\}/\mathbb{C}^{*}],$$ where the action is given by
$\lambda(x_0,\cdots,x_n)=(\lambda^{q_{0}}x_0,\cdots,\lambda^{q_{n}}x_n)$.
Thus we have proved the following proposition.
\begin{prop}
Every weighted projective stack of type $Q$ is a toric
Deligne-Mumford stack $\mathcal{X}(\mathbf{\Sigma})$ with an
underlying stacky fan $\mathbf{\Sigma}$. $\square$
\end{prop}

\begin{example}\label{example}
For $Q=(4,6,8)$, we have $d=2$ and $Q_{red}=(2,3,4)$. The weighted
projective stack $\mathbf{P}(Q_{red})={\bf P}_{2,3,4}$ has a
simplicial fan $\Sigma$ generated by $v_{0}=(-3,-2), v_{1}=(2,0),
v_{2}=(0,1)$, see \cite{Jiang}. We have the stacky fan
$$\mathbf{\Sigma}=(N,\Sigma,\beta),$$
where $N=\mathbb{Z}^{2}\oplus\mathbb{Z}_{2}$, and $\beta:
\mathbb{Z}^{3}\rightarrow N$ is given by the vectors
$$\{(2,0,1),
(0,1,0)),(-3,-2,0)\}.$$ The Gale dual map
$$\beta^{\vee}: \mathbb{Z}^{3}\rightarrow \mathbb{Z}$$
is given by the matrix $[4,6,8]$. So taking
$\mbox{Hom}(-,\mathbb{C}^{*})$ functor we have that the map
$$\alpha: \mathbb{C}^{*}\rightarrow (\mathbb{C}^{*})^{3}$$
is given by the matrix
$$\left[\begin{array}{c}
\lambda^{4}\\
\lambda^{6}\\
\lambda^{8} \end{array}\right].
$$
The toric Deligne-Mumford stack
$\mathcal{X}(\mathbf{\Sigma})=[\mathbb{C}^{3}\setminus
\{0\}/\mathbb{C}^{*}]$ is the weighted projective stack
$\mathbf{P}(Q)$ which is the nontrivial $\mu_{2}$-gerbe over
$\mathbf{P}(Q_{red})$.
\end{example}

\subsection{The weighted projective linear $2$-groups}\label{weighted2group}

From Section \ref{weightedstack} let $Z=\mathbb{C}^{n+1}\setminus
\{0\}$, then $\mbox{Aut}(Z)=GL(n,\mathbb{C})$, the nonsingular
$n\times n$ complex matrices over $\mathbb{C}$. Since the group
$\mathbb{C}^{*}$ embeds into the automorphism group
$\mbox{Aut}(Z)=GL(n,\mathbb{C})$ as a diagonal subgroup according
to the map $\alpha$ in the exact sequence (\ref{exactweighted3}),
we have the following Lemma.

\begin{lem}\Label{NorCen}
The centralizer $\widetilde{\mbox{Aut}^{0}(Z)}$ of $\mathbb{C}^{*}$
in  $\mbox{Aut}(Z)$ coincides with the normalizer
$\widetilde{\mbox{Aut}(Z)}$ of $\mathbb{C}^{*}$ in  $\mbox{Aut}(Z)$.
$\square$
\end{lem}

So we have a natural map $\varphi: \mathbb{C}^{*}\rightarrow \widetilde{\mbox{Aut}(Z)}$
and from the Introduction we have a $2$-group:
$$\mathbf{PGL}=[\mathbb{C}^{*}\stackrel{\varphi}\rightarrow \widetilde{\mbox{Aut}(Z)}]$$
which is called the weighted projective linear $2$-group. From the main result
of the paper or the Proposition in \cite{BN}, let $\mbox{Aut}(\mathbf{P}(Q))$
be the automorphism group of the weighted projective stack $\mathbf{P}(Q)$, we have:

\begin{cor}
The natural map $f: \mathbf{PGL}\rightarrow \mbox{Aut}(\mathbf{P}(Q))$ is
an isomorphism. $\square$
\end{cor}

\begin{example}
Let $B=\mathbf{P}^{d}$ be the $d$-dimensional projective space. We
give stacky fan $\mathbf{\Sigma}=(N,\Sigma,\beta)$ as follows. Let
$N=\mathbb{Z}^{d}\oplus \mathbb{Z}_{r}$ and $\beta:
\mathbb{Z}^{d+1}\longrightarrow N$ be the map determined by the
vectors:
$$\{(1,0,\ldots,0,0), (0,1,\ldots,0,0),\ldots,
(0,0,\ldots,1,0), (-1,-1,\ldots,-1,1)\}.$$
Then
$DG(\beta)=\mathbb{Z}$, and the Gale dual $\beta^{\vee}$ is given by
the matrix
$$[r,r,\ldots,r].$$
So we have the following exact
sequences:
$$0\longrightarrow \mathbb{Z}\longrightarrow \mathbb{Z}^{d+1}\stackrel{\beta}
{\longrightarrow}\mathbb{Z}^{d}\oplus\mathbb{Z}_{r}\longrightarrow
0\longrightarrow 0,$$
$$0\longrightarrow \mathbb{Z}^{d}\longrightarrow \mathbb{Z}^{d+1}\stackrel
{\beta^{\vee}}{\longrightarrow}\mathbb{Z}\longrightarrow
\mathbb{Z}_{r}\longrightarrow 0.$$ Then we obtain the exact
sequence:
\begin{equation}\Label{weightedsequence}
1\longrightarrow \mu_{r}\longrightarrow
\mathbb{C}^{*}\stackrel{\alpha}
{\longrightarrow}(\mathbb{C}^{*})^{d+1}\longrightarrow
(\mathbb{C}^{*})^{d}\longrightarrow 1.
\end{equation}
The toric Deligne-Mumford stack
$\mathcal{X}(\mathbf{\Sigma}):=[\mathbb{C}^{d+1}-\{0\}/\mathbb{C}^{*}]$
is the canonical $\mu_{r}$-gerbe over the projective space
$\mathbf{P}^{d}$ coming from the canonical line bundle, where the
$\mathbb{C}^{*}$ action is given by
$$\lambda\cdot
(z_{1},\ldots,z_{d+1})=(\lambda^{r}\cdot
z_{1},\ldots,\lambda^{r}\cdot z_{d+1}).$$ Denote this toric
Deligne-Mumford stack by $\mathcal{G}_{r}=\mathbf{P}(r,\ldots,r)$.
Let $Z=\mathbb{C}^{d+1}\setminus \{0\}$, we have that
$\mbox{Aut}(Z)=GL(d+1,\mathbb{C})$. The map $\alpha$ in
(\ref{weightedsequence}) is given by $\lambda\longmapsto
diag(\lambda^{r},\cdots,\lambda^{r})$, so
$\widetilde{\mbox{Aut}(Z)}\cong GL(d+1,\mathbb{C})$. The
automorphism 2-group of $\mathcal{G}_{r}$ is isomorphic to
$$\GG=[\mathbb{C}^{*}\stackrel{\varphi}{\longrightarrow}GL(d+1,\mathbb{C})]$$
where $\varphi$ is given by
$$\lambda\longmapsto\left[
\begin{array}{c}
\lambda^{r}\\
\vdots\\
\lambda^{r} \end{array} \right].
$$
\end{example}

\begin{example}
We give stacky fan $\mathbf{\Sigma}=(N,\Sigma,\beta)$ as follows.
Let $N=\mathbb{Z}\oplus \mathbb{Z}_{2}$ and $\beta:
\mathbb{Z}^{2}\longrightarrow N$ be the map determined by the
vectors:
$$\{(2,0), (-3,1)\}.$$
Then $DG(\beta)=\mathbb{Z}$, and the Gale dual $\beta^{\vee}$ is
given by the matrix $[6,4]$. So we have the following exact
sequences:
$$0\longrightarrow \mathbb{Z}\longrightarrow \mathbb{Z}^{2}\stackrel{\beta}
{\longrightarrow}\mathbb{Z}\oplus\mathbb{Z}_{2}\longrightarrow
0\longrightarrow 0,$$
$$0\longrightarrow \mathbb{Z}\longrightarrow \mathbb{Z}^{2}\stackrel
{\beta^{\vee}}{\longrightarrow}\mathbb{Z}\longrightarrow
\mathbb{Z}_{2}\longrightarrow 0.$$ Then we obtain the exact
sequence:
\begin{equation}\Label{weightedsequence2}
1\longrightarrow \mu_{2}\longrightarrow
\mathbb{C}^{*}\stackrel{\alpha}
{\longrightarrow}(\mathbb{C}^{*})^{2}\longrightarrow
(\mathbb{C}^{*})\longrightarrow 1.
\end{equation}
The toric Deligne-Mumford stack
$\mathcal{X}(\mathbf{\Sigma}):=[\mathbb{C}^{2}-\{0\}/\mathbb{C}^{*}]$
is the canonical $\mu_{2}$-gerbe over the weighted projective stack
$\mathbf{P}_{2,3}$ coming from the canonical line bundle, where the
$\mathbb{C}^{*}$ action is given by
$$\lambda\cdot
(z_{1},z_{2})=(\lambda^{6}\cdot z_{1},\lambda^{4}\cdot z_{2}).$$ The
weighted projective stack $\mathbf{P}_{4,6}$ is the moduli stack
$\overline{\mathcal{M}}_{1,1}$ of elliptic curves with one marked
point.  Let $Z=\mathbb{C}^{2}\setminus \{0\}$, we have that
$\mbox{Aut}(Z)=GL(2,\mathbb{C})$. The map $\alpha$ in
(\ref{weightedsequence2}) is given by $\lambda\longmapsto
diag(\lambda^{6},\lambda^{4})$, so $\widetilde{\mbox{Aut}(Z)}\cong
\mathbb{C}^{*}\times\mathbb{C}^{*}$. The automorphism 2-group of
$\mathbf{P}_{4,6}$ is isomorphic to
$$\GG=[\mathbb{C}^{*}\stackrel{\varphi}{\longrightarrow}\mathbb{C}^{*}\times\mathbb{C}^{*}]$$
where $\varphi$ is given by
$$\lambda\longmapsto\left[
\begin{array}{c}
\lambda^{6}\\
\lambda^{4} \end{array} \right].
$$
\end{example}


\end{document}